\documentclass{elsart}
\begin{document}
%\begin{center}

\begin{frontmatter}

\title{Countability of the Real Numbers}
%\thanksref{label1}}
\author[Slavica]{Slavica Vlahovic} 
%\address{Gunduliceva 2 Sisak, Croatia}
%\corauthref{cor1}}\thanksref{label2}}
and 
\author[Branislav]{Branislav Vlahovic\corauthref{cor}}
\ead{vlahovic@wpo.nccu.edu}
\corauth[cor]{Corresponding author.}

\address[Slavica]{Gunduliceva 2, Sisak, Croatia}
\address[Branislav]{North Carolina Central University, Durham, NC 27707, USA}

%\title{\bf Countability of the Real Numbers}

%\vspace {0.2truecm}

%\author{{Slavica Vlahovic$^1$ and Branislav Vlahovic$^{1,2}$ }\\
%{\normalsize $^1$Gunduliceva 2 Sisak, Croatia, $^2$North Carolina 
%Central University, Durham USA}}  
%\vspace {0.2truecm}

%\date{2 March 2004}
%\maketitle

\begin{abstract}

%\end{center}
%\vspace {0.2truecm}

%\hskip 2pc

The proofs that the real numbers are denumerable will be shown, i.e., that 
there exists one-to-one correspondence between the natural numbers $N$ and the 
real numbers $\Re$. The general element of the sequence that contains all real 
numbers will be explicitly specified, and the first few elements of the sequence will be written. Remarks on the Cantor's nondenumerability proofs of 1873 and 1891 that the real numbers are noncountable will be given.  

\end{abstract}

\begin{keyword}
% keywords here, in the form: 
denumerability \sep real  numbers \sep countability \sep cardinal numbers 


{\it MSC:} 11B05
\end{keyword}

\end{frontmatter}

\section{Introduction}

	The first proof that it is impossible to establish a one-to-one 
correspondence between the natural numbers $N$ and the real numbers $\Re$ is 
older than a century. In December 1873 Cantor first proved non-denumerability of 
continuum and that first proof proceeded as follows\cite{1,2,3,4}: Find a closed 
interval $I_0$ that fails to contain $r_0$ then find a closed subinterval $I_1$ 
of $I_0$ such that $I_1$ misses $r_1$ continue in this manner, obtaining an 
infinite nested sequence of closed intervals, $I_0 \supseteq I_1 \supseteq I_2 
\supseteq ...$, that eventually excludes every one of the $r_n$; now let $d$ be 
a point lying in the intersection of all the Ia's; $d$ is a real number 
different from all of the $r_n$. 

This proof that no denumerable sequence of elements of an interval [a,b] can 
contain all elements of [a,b] often is overlooked in favor of the 1891 diagonal 
argument\cite{5}, when reference is made to Cantor's proving the 
nondenumerability of the continuum. Cantor himself repeated this proof with some 
modifications\cite{2,3,6,7,8,9,10,11,12,13,14} from 1874 to 1897, and today 
we have even more variations of this proof given by other authors. However, we 
have to note that they are in nuce similar; all of them include same modification 
of the Cantor's idea to derive a contradiction by defining in terms which cannot 
possibly be in the assumed denumerable sequence. So, in principle, all these proofs 
do not represent a significant change from Cantor's original idea and we can take them to be the 
same as the Cantor's proofs.

	For the reason of clarity, we will not discuss objections to these 
proofs that have been raised earlier\cite{15,16,17,18,19,20,21} or the 
legitimacy of these proofs from intuitionistic points of view \cite{22} 
and their nonconstructive parts, namely appeal to the Bolzano-Weierstrass  
theorem\cite{23} and inclusion of impredicative methods\cite{24}.  We will focus 
to show what is in principle wrong with the general idea of Cantor's proofs and 
consequently all other proofs related to the statement that the set of real 
numbers is not denumerable.

	The main part of the paper is devoted to show that the real 
numbers are denumerable. The explicit denumerable sequence that contains all real numbers will be given. The general element that generates the sequence will be written as well as the first a few elements of that sequence. That there is 
one-to-one correspondence between the real numbers and the elements of the explicitly written sequence will be proven by the three independent proofs.  
	
	\section{Profs of the denumerability of the real numbers}

\vskip 0.5cm

{\bf Theorem 1}

\vskip 0.3cm

The real numbers $\Re$ are denumerable; it is possible to establish a one-to-one 
correspondence between the natural numbers $N$ and the real numbers $\Re$. In other words, the 
cardinal number $c$ of the set of real numbers is equal to the cardinal number  $\aleph_0$ of the set of natural numbers. The general element of the sequence that generates all elements of the set $\Re$ is as follows:
\begin{equation}
{a_1}^{{a_2}^{{a_3}^{.^{.^{.^{a_n}}}}}}
\label{gsequence}
\end{equation}
where in (\ref{gsequence}) each element $a_i$ of bases and exponents has the 
following form:
\begin{equation}
a_i={({m_{i1}\over n_{i1}})}^{[({{m_{i2}\over n_{i2}})}^{({m_{i3}\over 
n_{i3}})}]}
\label{gelement}
\end{equation}
where $m_{ij},n_{ij}\in N, i=1,2,3,...n, j = 1,2,3$. 

With the general element (1) it is possible to express each of the real numbers, and to generate the sequence which contains all real numbers. 
That can be done by writing (1) for all possible combinations of arguments, with the sum of all bases and exponents equal to 2,3, 4,... and so on. To do that in 
a way that will provide a one-to one correspondence between such produced set and the set of natural numbers $N$ all elements obtained by (1) can be for example arranged in the following way: a) For the fixed sum of bases and exponents write all possible fractions ${m_{11}\over n_{11}}$ of $a_1$ with lower denominator coming first. b) After that for the same fixed sum as before, write all elements $a_1={{{({m_{11}\over n_{11}})}^{({m_{12}\over n_{12}})}}},$ if it is possible to create such elements for that sum. Doing that, write first all possible combinations of ${m_{11}\over n_{11}}$  for fixed  ${m_{12}\over 
n_{12}}$  and only after that change ${m_{12}\over n_{12}}$  if it exists for that sum.  During that elements with lower denominator $n_{11}$ and $n_{12}$ will be written first again. c) If it is possible for that specific fixed sum of bases and exponents, following the same rules a) and b), continue by writing elements of shape 
$a_i={({m_{i1}\over n_{i1}})}^{[({{m_{i2}\over n_{i2}})}^{({m_{i3}\over 
n_{i3}})}]}$.
Again first change will be done in  ${m_{11}\over{n_{11}}}$ after that change of 
${m_{12}\over{n_{12}}}$ and at lastly the change of ${m_{13}\over{n_{13}}}$.  d) 
When all possible combinations of $a_1$ are written for the fixed sum of bases 
and exponents, continue with increasing the number of exponents, if it is 
possible for that sum, and continue by writing all combinations that correspond 
to ${a_1}^{a_2},  {a_1}^{{a_2}^{{a_3}}}, {a_1}^{{a_2}^{{a_3}^{{a_4}}}},...$ and 
so on.   During that, first change exponents $a_i$ with lower index.
e) When all possible increases of exponents are done and all possible 
combinations for the fixed sum are written, increase the value of the sum and 
repeat procedures a) through e). 
In addition all elements that appear again will not be written down. First few 
elements of this sequence are as follows:
$${1\over1},  {2\over1},  {1\over2}, {3\over1}, {1\over3}, {4\over1},  
{3\over2},  
{2\over3}, {1\over4},  {5\over1},  {1\over5},  {({2\over1})}^{1\over2},  
{({1\over2})}^{1\over2}, 
{6\over1},  {5\over2},  {4\over3},  {3\over4},   {2\over5},  {1\over6},  
{({3\over1})}^{2\over1},  {({1\over3})}^{2\over1},   {({3\over1})}^{1\over2},  
{({1\over3})}^{1\over2}, $$ 

$${({2\over1})}^{3\over1},  {({1\over2})}^{3\over1},  {({2\over1})}^{1\over3}, 
{({1\over2})}^{1\over3},   {7\over1},   {5\over3},  {3\over5},   {1\over7},  
{({4\over1})}^{2\over1},  {({3\over2})}^{2\over1}, {({2\over3})}^{2\over1},     
{({1\over4})}^{2\over1},  {({3\over2})}^{1\over2},  {({2\over3})}^{1\over2},
{({3\over1})}^{3\over1},  {({1\over3})}^{3\over1},  {({3\over1})}^{1\over3},  
{({1\over3})}^{1\over3},$$

$${({2\over1})}^{3\over2},  {({1\over2})}^{3\over2},  {({2\over1})}^{1\over4},  
{({1\over2})}^{1\over4},  {7\over2},  {5\over4},  {4\over5},  {2\over7},  
{({5\over1})}^{2\over1}, {({1\over5})}^{2\over1},  {({5\over1})}^{1\over2},  
{({1\over5})}^{1\over2},  {({4\over1})}^{3\over1},  {({3\over2})}^{3\over1},    
{({2\over3})}^{3\over1},  {({1\over4})}^{3\over1},  {({4\over1})}^{1\over3},$$

$${({3\over2})}^{1\over3},  {({2\over3})}^{1\over3},  {({1\over4})}^{1\over3},  
{({3\over1})}^{4\over1},  {({1\over3})}^{4\over1},  {({3\over1})}^{3\over2}, 
{({1\over3})}^{3\over2},  {({3\over1})}^{2\over3},  {({1\over3})}^{2\over3}, 
{({3\over1})}^{1\over4},  {({1\over3})}^{1\over4}, {({2\over1})}^{5\over1},  
{({1\over2})}^{5\over1},  {({2\over1})}^{1\over5},  {({1\over2})}^{1\over5},$$
 
 \begin{equation}
 {({2\over1})}^{[({2\over1})^{1\over2}]},   
{({1\over2})}^{[({2\over1})^{1\over2}]},   
{({2\over1})}^{[({1\over2})^{1\over2}]}, 
{({1\over2})}^{[({1\over2})^{1\over2}]},...      
\label{gexelement}
\end{equation}

\vskip 0.5cm
{\bf Proof of the theorem 1}

\vskip 0.3cm

Note first that it is obvious that the sequence contains all rational and 
algebraic numbers, and that transcendental numbers are included also, as in the case 
of algebraic irrational exponents and algebraic bases [25-31], for instance, for                                
    \begin{equation}
    {({2\over1})}^{[({2\over1})^{1\over2}]}=2^{\sqrt{2}}
    \label{geselement}
    \end{equation}
 
Exponents $a_i$ in the general element (1) can be either algebraic or 
transcendental, which depends on arguments $m_{i,j},n_{i,j}$ of $a_i$. How
$a_i$ has the shape (2) and the arguments of $a_i$ can be changed for an arbitrary small amount, it is obvious that $a_i$ can obtain a value in any chosen interval.  Since the general elements of the sequence (1) have the form 
${a_1}^{{a_2}^{{a_3}^{.^{.^{.^{a_n}}}}}},$ 
hence exponential function is continuous, and because values of the arguments of (1), $a_1, a_2, a_3,..., a_n$ can be chosen from any interval and can be changed 
independently one from another for an arbitrary small amount, it follows that 
expression (1) can obtain any arbitrary value.  Therefore, with (1), in any arbitrary chosen interval one can generate infinitely many algebraic and transcendental numbers, which is actually the continuum [32].  That means that with (1) we can represent any real number, which proves the theorem 1. 

Since this is an extremely important issue, we will give two additional 
completely independent proofs of the theorem. However, before that let us 
consider some properties of the sequence generated with (1). 

We need to note that only with $a_1$, which also can be transcendental, such as 
in (\ref{geselement}),
it is not possible to express all numbers, for instance the number e, as it 
requires
\begin{equation}
e={m_1\over n_1}^{[{m_2\over n_2}^{m_3\over n_3}]}
 \label{e}
    \end{equation}
that is
\begin{equation}
1={m_2\over n_2}^{m_3 \over n_3}ln{m_1\over n_1}
 \label{le}
    \end{equation}
and this cannot be, because $ln{m_1\over n_1}$ is always transcendental [31-33] for $m_1,n_1\in N$. However, it is necessary to note that no reason exists that could prevent expressing any arbitrary number with ${a_1}^{a_2}$.   
Therefore it may already be possible by ${a_1}^{a_2}$ to express all real numbers, i.e. maybe it is not necessary to build numbers with more and more exponents, i.e. numbers of the shape
       \begin{equation}
{a_1}^{{a_2}^{{a_3}}}, {a_1}^{{a_2}^{{a_3}^{a_4}}} ,...,
{a_1}^{{a_2}^{{a_3}^{.^{.^{.^{a_n}}}}}},...
 \label{geseries}
    \end{equation}
However, for now that statement can not be established, because 
it is not possible for now to calculate ${a_1}^{a_2}$ for the general case when both $a_1 $ and $a_2$ are transcendental [34], it is not even possible to calculate it for earlier simpler case (4). 

We need also to note the following: the statement that ${a_1}^{a_2}$ has only 
$\aleph_0\aleph_0 = \aleph_0$ elements and that this is the reason why it cannot 
contain all real numbers, which we have ${\aleph_0}^{\aleph_0}$ is not a good 
argument, because $\aleph_0= {\aleph_0}^{\aleph_0}$ if the set of real number 
is countable, as it is. Consequently the possibility that the set of all 
real numbers could be expressed by only ${a_1}^{a_2}$ must be kept open.

\vskip 0.3cm
{\bf Theorem 2}
\vskip 0.2cm

	The set of numbers, generated by general element (1) and procedure a) 
through e) given in theorem (1), does not have any gaps. At each cut of the set 
the first component of the cut has the last element, or the second component of 
the cut has the first element, or both of these cases occur.

\vskip 0.2cm
{\bf The proof of the theorem 2}
\vskip 0.2cm

	In the definition of the theorem the meaning of the cut is simply a rule 
for dividing a set into two non-empty parts $A$ and $B$ such that every element 
of $A$ precedes every element of $B$ while $A$ and $B$ together exhaust the set.

	Let denote with $S$ the set generated by (1) followimng the procedure given in 
theorem 1. The set $S$ obviously has as a subset the set of algebraic numbers 
$R$. If the $A/B$ is the cut in $S$, where $A$ is the first component and $B$ 
the second component of the cut, then $(A\cap R)/(B\cap R)$ is the certain
defined cut $k$. If the component $A$ of the cut $A/B$ has the last 
element, then the cut does not generate new elements, $k \in A$ and $k$ is the 
last element of $A$. 

	If the $k$ is not the last element in $A$ then exists 
\begin{equation}
k^\prime \in A \quad such \quad that \quad k<k^\prime.
\label{kprime}
    \end{equation}
 But it is not possible that $k^\prime \in R$, since it 
requires that $k^\prime 
\in B\cap R$ and consequently $k^\prime \in B$, which is not possible because of 
(\ref{kprime}), since $A$ and $B$, as components of the cut do not have the 
common elements. So, if $k^\prime 
\not\in R$, than it means that $k^\prime$ is a given cut $C/D$, of the set R, 
created 
with a gap in the set $R$ and for that reason the first component $C$ of the cut 
$C/D$ does not have the greatest element. Since $k<k^\prime$, it follows that 
$A\cap R 
\subset C$. If $k^{\prime\prime}$ is any element of the set $C\backslash A$ then 
$k<k^{\prime\prime}$ will require that $k^{\prime\prime}\in B$. But because of 
the $k^{\prime\prime}\leq k^\prime$, $k^\prime\in A$ it 
requires $k^{\prime\prime}\in A$. Both relations $k^{\prime\prime}\in A$ and 
$k^{\prime\prime}\in B$ can not be 
satisfied, because $A$ and $B$ as the components of the cut in $S$ are 
disjunctive 
sets. So, if $k\in A$ then $k$ is the last element of the component $A$.  

	In the same way it can be proven that if the $k\in S\backslash A$, so 
$k\in B$, then the $k$ is the first element of the second component $B$. 

	This proves the theorem. It is important to note that in proving this 
theorem we used the following properties of the set $S$: a) that it has the 
dense subset of the algebraic numbers $R$ and b) the set $S$ is everywhere dense, 
consequently for any $C\backslash A$ the general element of the set, 
relation (1), will generate numbers $k^{\prime\prime}\in 
C\backslash A$, which are required to be in both $A$ and $B$, which 
established the contradiction. 
	
It is not necessary to note that this also proves that the set $S$ is equivalent 
to the set of all real numbers $\Re$, since the set does not have the first and 
the last element, it is dense, it does not have any gaps, and it is linear. 
	
 \vskip 0.3cm
{\bf Theorem 3}
\vskip 0.2cm
 
	The set $S$, generated by the general element (1) and procedure a) 
through e) given in theorem 1, is similar (isomorphic) to the set of all real 
numbers $\Re$. 

\vskip 0.2cm
{\bf The proof of the theorem 3}
\vskip 0.2cm

The set $\Re=(\Re:<)$ of all real numbers ordered by the magnitude of its 
elements has the following properties:
a)	it does not have the first and the last element,
b)	it is continuous in Dedekin's sense, and
c)	it is separable.
Any other set with properties a) to c) is similar to the set of real numbers 
$\Re$. 
Let us now show that the set $S$ given with theorem 1 and general 
element 
(1) satisfies the conditions a) to c) and is similar to the set of real numbers 
$\Re$.  The set $S$ obviously satisfies properties a) and c). It also satisfies 
property given by b), as it is proven by theorem 2. With that the theorem 3 is 
proven. 
However, to keep this proof independent from the theorem 2, we will now prove it 
without using that theorem.

Let us denote by $M_1 \subseteq \Re$ any countable part of set $\Re$, such 
that it satisfies condition: d) each interval of $\Re$ contains at least 
one element of $M_1$ and each interval of $M_1$ contains at least one 
element from $\Re\backslash M_1$. It is obvious that the set $M_1$ with the property d) can be 
created, since between any two algebraic numbers exist at least one transcedental
number and between any two transcedental numbers exist at least one algebraic number 
[35].  
	
	Let us denote by $M_2 \subseteq S$ any countable part of $S$, such that 
each interval of $S$ contains at least one element of $M_2$ and each interval of 
$M_2$ contains at least one element of $S\backslash M_2$.  The set $M_2$ also can 
be created, since the set $S$ generated by (1) also obviously has the property 
d). The set $S$ has as a subset the set of algebraic numbers.
% $m_{ij}\over n_{ij}$. 
 Also between any arbitrary chosen pairs of algebraic numbers it is possible 
to create transcedental numbers by the general element of sequence (1), and 
between any pairs of transcedental numbers generated by (1) there are algebraic 
numbers generated by (1).
% $m_{ij}\over n_{ij}$. 

	Further, the sets $M_1$ and $M_2$ satisfy the following: a) sets do not have the first 
or the last element, b) sets are dense, c) sets are countable. 

The sets $M_1$ and $M_2$ are similar to the set of algebraic numbers. We will 
now show that any similarity 
\begin{equation}
\varphi (x), (x\in M_1), \varphi (M_1) = M_2 
\label{similarity}
\end{equation}
between $M_1$ and $M_2$ can be extended on the similarity between entire $\Re$ and 
$S$.

Let take $x\in \Re\backslash M_1$ then we have cut 
\begin{equation}
M_1 = (-\infty,x)_{M_1}\cup(x,\infty)_{M_1} 
\label{cutm1}
\end{equation}
in the set $M_1$, which because of the density of the set $M_1$ opens 
a gap in $M_1$ and an element $x\in \Re\backslash M_1$ fulfills that gap in 
$\Re$. The $x$ is actually the only element that is between summands 
(\ref{cutm1}).

The cut in the set $M_1$ by the similarity (\ref{similarity}) makes cut 
 \begin{equation}
M_2 = \varphi (M_1) = \varphi (-\infty,x)_{M_1}\cup \varphi (x,\infty)_{M_1} 
\label{cutm2}
\end{equation}
 of the set $M_2$. Because of the similarity of the sets $M_1$ and $M_2$, the 
cut (\ref{cutm2}) creates the gap in $M_2$. In that gap, because of the property 
d) is the element of the set $S$, which is defined by the similarity between 
$M_1$ and $M_2$, and by the element $x\in \Re\backslash M_1$, i.e. by the 
$\varphi (x)$.  By that the transformation $\varphi (x)$  is defined for each $x\in \Re$. 
Obviously $\varphi (\Re) = S$. With that the theorem is proven, the set 
$S$ is similar to the set $\Re$.

\section{Remarks on the Cantor's proofs}

The above proposed sequence that contains all real numbers, and established denumerability of the real numbers are obviously in contradiction with the Cantor's proofs of nondenumerability. It is not to us to find the errors in Cantor's proofs and all numerous variations of his proofs that currently exist. However, we will give remarks on the Cantor's two most quoted proofs, from 1873 and 1891.     
         
 \vskip 0.3cm
 {\bf Theorem 4}
	\vskip 0.2cm
	
	In the Cantor's 1873 proof of nondenumerability, Cantor stated that it is possible to create 
sequences of progression and 
	regression of elements, which allow for any interval of real numbers 
$(\alpha...\beta)$ to define, in the limit, a number $\eta\in(\alpha,\beta)$, 
which was not included in the sequence assumed to contains all real numbers. The existence of the limit $\eta$ does not lead to the conclusion that the number $\eta$ is not in the sequence assumed to countain all real numbers and that the set of all real numbers is not countable.
	
	\vskip 0.2cm
	{\bf Proof of the theorem 4}
	\vskip 0.2cm
	
	Let us now look in Cantor's 1873 nondenumerability proof, which appeared 
in Crelle's Journal in January 1874. 
	
	Assuming that the real numbers are countable, it follows that they could 
be sequenced on an index of natural number $N$: 
	\begin{equation}
	  \omega_1, \omega_2, \omega_3, ..., \omega_\nu,...  
	  \label{om1}
\end{equation}

	  Cantor then stated that for any given interval $(\alpha...\beta)$ he 
could show the existence of a number $\eta\in(\alpha,\beta)$ which is not 
included in the sequence (\ref{om1}).
	  
	  Assuming $\alpha < \beta$, he picked the first two numbers from 
(\ref{om1}), which fell within the interval $(\alpha,\beta)$. Denoted 
$\alpha',\beta'$,  respectively, these were used to constitute another interval 
$(\alpha'...\beta')$. Proceeding analogously, Cantor provided a sequence of 
nested intervals, reaching $(\alpha^{(\nu)}...\beta^{(\nu)})$, where 
$\alpha^{(\nu)}, \beta^{(\nu)}$ were the first two numbers from  (\ref{om1}) 
lying within $(\alpha^{(\nu-1)}...\beta^{(\nu-1)})$
	  
	  If the number of intervals thus constructed were finite, then at most 
only one more element from (\ref{om1}) could lie in 
$(\alpha^{(\nu)},\beta^{(\nu)})$. It was easy in this case for Cantor to 
conclude that a number $\eta$ could be taken in this interval which was not 
listed in  (\ref{om1}). Clearly any real number $\eta\in 
(\alpha^{(\nu)},\beta^{(\nu)})$ would suffice, as long as $\eta$ was not the one 
element possible listed in (\ref{om1}). 
	  
	  In the case when the number of intervals 
$(\alpha^{(\nu)},\beta^{(\nu)})$ were not finite, Cantor's argument shifted to 
consider two alternatives in the limit. Since the progressing sequence $\alpha, 
\alpha',...,\alpha^{(\nu)},...$ did not increase indefinitely, but was bounded 
within $(\alpha,\beta)$, it had to assume an upper limit which Cantor denoted 
$\alpha^\infty$. Similarly, the regression sequence $\beta, 
\beta',...,\beta^{(\nu)},...$ was assigned the lower limit $\beta^\infty$. Where 
$\alpha^\infty < \beta^\infty$, then, as in the finite case, any real number 
$\eta\in(\alpha^\infty,\beta^\infty)$ was sufficient to produce the necessary 
real number not listed in  (\ref{om1}).  However, were 
$\alpha^\infty=\beta^\infty$, Cantor reasoned that 
$\eta=\alpha^\infty=\beta^\infty$ could not be included as an element of 
(\ref{om1}) (we will prove that this assumption is not correct). 
He designed $\eta=\omega_\rho$. But $\omega_\rho$, for sufficiently 
large index $\nu$, would be excluded from all intervals nested within 
$(\alpha^{(\nu)},\beta^{(\nu)})$. Nevertheless, by virtue of the construction 
Cantor had given, $\eta$ had to lie in every interval 
$(\alpha^{(\nu)},\beta^{(\nu)})$, regardless of index. The contradiction 
established the proof: $R$ was nondenumerable. 
	  
	  The main part of the proof is that there is a progression of elements 
$\alpha^{(n)}$ and regression of elements $\beta^{(n)}$, such that 
	
	   	\begin{equation}  
	\alpha<\alpha^{(1)}<\alpha^{(2)}<...<...<\beta^{(2)}<\beta^{(1)}<\beta
	  \label{om2}
\end{equation}

The progression ought to have an upper limit; but there is no element 
$\alpha^{(n)}$ which can serve as this upper limit, for if any element 
$\alpha^{(n)}$ is proposed, one can clearly carry the process just indicated 
that $\alpha^{(n)}$ will be outside the interval $\alpha^{(n)}...\beta^{(n)}$.

The best way to illustrate what is wrong with this proof is to apply it on the 
set of all rational numbers. Applying exactly the same procedure proposed by 
Cantor on the set of rational numbers from interval (0,2) it is possible to make 
the sequences that determine progression of elements $\alpha^{(n)}$, and 
regression of elements $\beta^{(n)}$ such that

\begin{equation}   
0<{2\over4}<{4\over6}<...<{{2+2n}\over{4+2n}}<...<...<{{4+2n}\over{2+2n}}<...<{8
\over6}<{6\over4}<2
   \label{om3}
\end{equation}

  The above progression and regression of elements determine as the limit number $\eta = 1$.  There is no element $\alpha^{(n)}$ from (\ref{om3}) which can serve as this upper limit. Following the Cantor's line of conclusion, for again 
$\eta=\alpha^\infty \in (\alpha^{(n)},\beta^{(n)})$ for all $n$, and hence 
$\alpha^\infty \neq \alpha_n$ for all n, simply by the way the progression and 
regression sequences are constructed the number $\eta=1$, which is obtained in 
the limit, cannot be in the sequence of rational numbers on interval  
$(0,2)$. Therefore we should conclude that the set of 
rational numbers is not denumerable, while we know that it is denumerable. Why we 
get this contradiction?  

%To describe that let look step by step in the Cantor's proof. 
%  
%  First Cantor assumed that a given continuous sequence to be denumerable. Then 
%without disturbing the order of the elements he attached to each element a 
%definite natural number $\alpha^{(n)}$. Then he assumed that without the loss of 
%generality the elements have been so  numerated that the element $\alpha^{(1)}$ 
%presides the element $\alpha^{(2)}$, etc. This is required to make regression 
%sequence. However, this assumption is not correct. The statement that numbers 
%can be denumerated does not imply that they may be denumerated by magnitude. For 
%instance we know how to make set of rational numbers countable, but we can not 
%denumerate them by their values, because the set of rational numbers is dense.   

Answer is simple. The above example demonstrates that the Cantor statement that the number 
$\eta$ which is obtained as the limit of progression and regression sequences cannot be an element 
of (\ref{om1}) is not correct. He stated that whatever $\omega_\rho$ is taken for $\eta$ that for 
sufficient large index $\nu$, it will be excluded from all intervals nested within 
$(\alpha^{(\nu)},\beta^{(\nu)})$. Obviously there is no element in sequences (\ref{om2}) or (\ref{om3})
which can serve as the limit $\eta$.  Any number $\alpha^{(n)}$ taken from (\ref{om2}) or (\ref{om3}) will fail 
as Cantor properly stated. However, it is not correct that there is no number from (\ref{om1}), which is equal to $\eta$ and which will be inside any interval $(\alpha^{(\nu)},\beta^{(\nu)})$ for any $\nu$.
In our example number $1$ is obviously in (\ref{om1}) and it is in all intervals of (\ref{om2}) or (\ref{om3}) regardless of how large is the $\nu$. So, creating nested intervals $(\alpha^{(\nu)},\beta^{(\nu)})$ by following Cantor's 
procedure, as the result of regressing and progressing sequences, obtained is in the limit number $\eta$, 
which by his statement cannot be part of (\ref{om1}), because of the way how it is created. This is obviously not correct, since as our example demonstrates, the number $\eta=1$ which is obtained following Cantor's procedure is 
obviously part of the sequence (\ref{om1}), since in our example (\ref{om1}) represents sequence of rational numbers and $\eta=1$ is the part of that sequence.        
By getting for the limit in sequence (\ref{om2}) a number which is not an element 
of (\ref{om2}) it does not mean that the set of all real numbers is not 
countable. It only means that the particular sequence does not contain that 
number, but the same number may be an element of the sequence (\ref{om1}). The 
particular sequence (\ref{om2}) is not the only sequence that can be constructed 
from (\ref{om1}), so it does not need to contain all numbers from (\ref{om1}).

 \vskip 0.3cm
{\bf Theorem 5}
\vskip 0.2cm

By the Cantor's diagonal procedure, it is not possible to build numbers that are 
different from all numbers in a general assumed denumerable sequence of all real 
numbers. The numbers created on the diagonal of the assumed sequence have the 
values that are not different from the values of the numbers in the assumed 
denumerable sequence. 

\vskip 0.2cm
	{\bf Proof of the theorem 5}
	\vskip 0.2cm
	
In his proof Cantor first produced a countable listing of elements $E_\nu$ in 
terms of the corresponding array (\ref{ar1}), where each $a_{\mu,\nu}$ was 
either m or w:

%\begin{center}
 \begin{equation}
 {\begin{array}{clcr}
     E_1=(a_{11},a_{12},...,a_{1\nu},...)\\
     E_2=(a_{21},a_{22},...,a_{2\nu},...)\\ 
     \vdots 
\;\;\;\;\;\;\;\;\;\;\;\;\;\;\;\;\;\;\;\;\;\;\;\;\;\;\;\;\;\;\;\;\;\;\; \\
     E_{\nu}=(a_{\mu1},a_{\mu2},...,a_{\mu\nu},...)\\
     \vdots 
\;\;\;\;\;\;\;\;\;\;\;\;\;\;\;\;\;\;\;\;\;\;\;\;\;\;\;\;\;\;\;\;\;\;\;
     \end{array}}
       \label{ar1}
 \end{equation}
 %\end{center}
       
       Then he defined a new sequence $b_1, b_2,...,b_\nu,...$ . Each $b_\nu$ 
was either m or w, determined so that $b_\nu\neq a_{\nu\nu}$. By formulating 
from this sequence of $b_\nu$ the element $E_0=(b_1,b_2,...,b_\nu,...)$, it 
followed that $E_0 \neq E_\nu$ for any value of the index $\nu$.  

Let apply Cantor's procedure on the set of real numbers from the interval 
(0,1) to answer what has in reality Cantor proved by his diagonal 
procedure. He claimed that it is possible on the diagonal of an arbitrary 
denumerable sequence, which represents numbers in the interval (0,1), to create 
numbers that are different from the first number in the first decimal point, 
that are different from the second number in the second decimal point, and so 
on. From that he concluded that the created numbers are not in that assumed 
sequence and that the real numbers are not denumerable. The idea of his proof is 
that in the arbitrary assumed denumerable set of real numbers, each element of 
the set, each real number, has to be related by a one-to one correspondence to a 
natural number, which has a final value. That is, any real number has to take a 
finite place in that sequence. After that he concluded that numbers on the 
diagonal will be different from any of numbers in the sequence, because they are 
different from the first, the second and other numbers. But has Cantor really 
proved this, and what does it means that a created number on the diagonal is 
different from the first number in the sequence, the second and so-on? By the 
way how a diagonal number is created it is obvious that it is different from the 
first number in the sequence in $10^{-1}$ of the magnitude, from the second on 
$10^{-2}$ and from an $n^{th}$ number on $10^{-n}$ order of the magnitude. Of 
course, it will be different from any $n$ elements in that sequence. From this and 
the earlier statement that any real number has to be assigned to a finite 
position in the sequence, Cantor concluded that the numbers created on the diagonal are different from all numbers in the arbitrary assumed sequence because any number has to be assigned to an $n$ that has a final value. What Cantor has proved by the diagonal procedure is that if we take any finite subset of the real set, i.e. first $n$ elements, it is possible to create an $n+1$ element by the defined diagonal procedure. It is true that with this procedure we can go further and build more and more elements of the real set, which are different from an finite subset that contains $n$ elements. However, he did not prove that it is possible to create new numbers that are not already included in the arbitrary assumed 
countable set, the numbers that will be different from all elements in that set. 
The conclusion that created diagonal elements are such numbers is not correct. 
A number created on the diagonal has to be different not only from the first 
$n$ elements in the set, regardless of how large is the $n$. The created number must 
be different from all numbers in that set, which as we know has an infinite number 
of elements. This is the main point, that the assumed denumerable sequence has an 
infinite, not finite, number of elements.  
%not only from an finite subset regardless how large that subset is.  
However, by the way how the diagonal numbers are created it is obvious that they 
are different just from final subset from the assumed denumerable set and not from all 
numbers in that set.  To prove that let us look what is the difference in the 
value between the numbers created on the diagonal and the numbers included in the 
assumed denumerable set. The difference between these diagonal numbers and all other 
numbers in the assumed sequence is simply given by the equation:

\begin{equation}
         \lim_{n \rightarrow \infty} 10^{-n} = 0         
  \label{lim2}
 \end{equation}
         
In the above equation we have to take the limit when $n$ is going to infinity, since 
we have to take into account that the proposed denumerable set has infinite 
number of elements. The numbers created on the diagonal must be different from 
all these numbers, not only from a final size subset of these elements. 

This proves that the numbers created on the diagonal are not different from all 
the numbers in that assumed denumerable set. They are different from a subset with 
$n$ elements, but they are contained in the proposed denumerable set. 
Cantor did not took into account that the proposed denumerable set has infinite 
number of elements. In his discussion he was always focused on some finite size subset, 
considering some finite number of elements that belong to some finite $n$.
It is impossible by the proposed diagonal procedure to build numbers that are not 
included in the assumed denumerable set and particularly it is not possible by this way 
to create an ascending hierarchy, in fact a limitless sequence of transfinite 
powers.

\section{Conclusion}

It is shown that the set of all real numbers is denumerable. The general element 
that generates the set is given and the first few elements of the sequence that 
contains all real numbers are written explicitly. By three independent proofs it 
is shown that the proposed sequence represents the set of the numbers which is 
dense anywhere, that does not have any gaps, and that is similar to the set of all real numbers, which proves that the sequence contains all real numbers. It is also proven that the Cantor's 1873 proof of non denumerability is not correct since it implicates non denumerability of rational numbers. In addition it is proven that the numbers generated by the diagonal procedure in Cantor's 1991 proof are not different from the numbers in the assumed denumerable set.

  	  	\end{document}